\documentclass{amsart}

\newtheorem{proposition}{Proposition}

\newtheorem{lemma}[proposition]{Lemma}
\newtheorem{corollary}[proposition]{Corollary}

\newtheorem{remark}[proposition]{Remark}

\begin{document}
\title{Generic fiber of  power series ring extensions}
\author{Tiberiu Dumitrescu}
\address{Facultatea de Matematica si Informatica, Universitatea Bucuresti,
Academiei Str. 14, Bucharest, RO-010014, Romania}
\email{tiberiu\underline{ }dumitrescu2003@yahoo.com, tiberiu@fmi.unibuc.ro}
\date{}
\maketitle
\begin{abstract}
Let $D$ be a Noetherian domain  containing a field,
$a\in D$  a nonzero nonunit and  $z$  an indeterminate over $D$.
We prove that the generic fiber of the
extension $D[[z]]\hookrightarrow D[1/a][[z]]$ has
dimension {\em $\geq \mbox{dim}(D/aD)$.}
\end{abstract}\vspace{4mm}
\vspace{10mm}\\
{\bf 1. Introduction.}\\

An extension $A\subseteq B$ of integral domains
is said to be a {\em trivial generic fiber} (TGF) extension, if
every nonzero prime ideal of $B$ has nonzero intersection with $A$.
Let $K$ be a field and $z,x_1,x_2,...,x_n$ indeterminates over $K$.

In \cite{HRW}, Heinzer, Rotthaus and Wiegand proved that
the mixed polynomial/po\-wer series ring extension
$$K[x_1,x_2,...,x_n][[z]]\hookrightarrow K[x_1,x_2,...,x_n,1/x_1][[z]]$$
is  TGF for $n=1$ and non-TGF for $n=3$.
They asked what happens for $n=2$.

In this note, we study the TGF property for power series extensions
of type $$D[[z]]\hookrightarrow D_a[[z]]$$ where
$D$ is a domain, $0\neq a\in D$ and $z$ an indeterminate over $D$.
Here  $D_a$  is a (short) notation for  fraction ring $D[1/a]$.
We prove that $D[[z]]\hookrightarrow D_a[[z]]$ is TGF when $D$ is one-dimensional.

Let $A$ be a domain
and $B$ a ring between $A[x,y]$ and $A[[x,y]]$, where $x,y,z$ are indeterminates over $A$.
We show that the extension $B[[z]]\hookrightarrow B_x[[z]]$ is not TGF.
In particular, for a field $K$, the extension
$$K[x,y][[z]]\hookrightarrow K[x,y]_x[[z]]$$
is not TGF, thus aswering \cite[Question 4.8]{HRW} mentioned above.
The key point of our proof is that,
as shown by Abhyankar and Moh in \cite[Theorem 2]{AM},
there exist  elements  $\lambda\in zK[x]_x[[z]]$ which are analytically independent over
$K[[x,z]]$ (e.g., $x\sum_{n\geq 1} (z/x)^{n!}$).
Then it is easy to find a prime ideal of $K[x,y]_x[[z]]$ containing $\lambda-y$
which has zero intersection with $K[x,y][[z]]$.

The same kind of arguments
can be used to show that the generic fiber of the extension
$$K[x,y_1,...,y_n][[z]]\hookrightarrow K[x,y_1,...,y_n]_x[[z]]$$
has dimension $\geq n$.

Finally, using this fact and Cohen's structure theorem for complete local rings,
we establish the following result.
Let $B$ be a Noetherian domain  containing a field,
$a\in B$  a nonzero nonunit and  $z$  an indeterminate over $B$.
Then the generic fiber of the extension $B[[z]]\hookrightarrow B_a[[z]]$ has
dimension {\em $\geq \mbox{dim}(B/aB)$.}
In particular, $B[[z]]\hookrightarrow B_a[[z]]$ is not TGF
when {\em $\mbox{dim}(B/aB) \geq 1$.}
\\[10mm]
{\bf 2. Results.}\\

Let $K$ be a field and $x,z$ indeterminates over $K$.
In \cite[Proposition 2.6]{HRW}, it is shown that
the  extension $K[x][[z]]\hookrightarrow K[x]_x[[z]]$ is  TGF.
We extend this result. 

\begin{proposition}\label{2}
Let $A$ be a one-dimensional domain, $a\in A$ a nonzero nonunit and $x$ an indeterminate over $A$.
Then the extension $A[[x]]\hookrightarrow A_a[[x]]$ is TGF.
\end{proposition}
{\noindent\it Proof.}
We adapt the proof of \cite[Proposition 2.6]{HRW}.
Suppose that the extension is not TGF.
Hence there exists a nonzero prime ideal $P$ of $A_a[[x]]$ such that $P\cap A[[x]]=0$.
Since $x\not\in P$, there exists an $f\in P$ such that $0\neq f(0)\in A$.

Let $u$ be the canonical map
$A[[x]]\rightarrow A_a[[x]]/fA_a[[x]]$. Reducing $u$ mo\-du\-lo $(x)$ we get the canonical map
$A\rightarrow A_a/f(0)A_a$ which is surjective because $A/f(0)A$ is zero-dimensional, cf. \cite[Theorem 3.1]{H}.
It is well-known that $A[[x]]$ is complete in the $(x)$-adic topology. Moreover,
by the next lemma, $A_a[[x]]/fA_a[[x]]$ is separated in the $(x)$-adic topology.
By Cohen's theorem (see \cite[Lemma on page 212]{M} or \cite[Lemma 4.2]{BRD}),
$u$ is surjective.
Since $fA_a[[x]]\cap A[[x]]\subseteq P\cap A[[x]]=0$, $u$ is an isomorphism.
This is a contradiction, because $a$ is invertible in $A_a[[x]]/fA_a[[x]]$ but non-invertible in  $A[[x]]$.
$\bullet$\\[2mm]

The following lemma is probably well-known, but we were unable to find a reference for it.

\begin{lemma}\label{1}
Let $A$ be a  domain and $0\neq f\in A[[x]]$.
Then the principal ideal $fA[[x]]$ is closed in the $(x)$-adic topology.
\end{lemma}
{\noindent\it Proof.}
Let $g\in \cap_{m\geq 1} (f,x^m)A[[x]]$.
Clearly, we may assume that $c=f(0)$ is nonzero. Then $g=fh$ for some $h\in A_c[[x]]$,
say, $h=\sum_{n\geq 0}h_nx^n$ with $h_n\in A_c$.
Let $k\geq 1$. As $g\in (f,x^k)A[[x]]$, $g=fq+x^kr$ for some $q,r\in A[[x]]$.
So $g=fh=fq+x^kr$, hence $f(h-q)=x^kr$. As $f(0)\neq 0$, $h-q=x^ks$ for some $s\in A_c[[x]]$.
Hence $h_0,h_1,...,h_{k-1}\in A$. So $h\in A[[x]]$, that is, $g\in fA[[x]]$.
$\bullet$\\[2mm]

Let $R\subseteq T$ be an extension of domains and $\lambda_1$,...,$\lambda_n\in zT[[z]]$,
where $z$ is an indeterminate over $T$.
Recall that $\lambda_1$,...,$\lambda_n$ are said to be {\em analytically independent} over $R[[z]]$ if the
$R[[z]]$-algebra morphism
$$\theta:R[[z]][[y_1,...,y_n]]\rightarrow T[[z]]$$
given by $\theta(y_i)=\lambda_i$, $1\leq i\leq n$, is injective,
where $y_1,...,y_n$ are indeterminates over $R[[z]]$.
The next proposition is the key technical result of this note.

\begin{proposition}\label{3}
Let $A\subseteq A'$ be an  extension of domains,  $a\in A$ a nonzero nonunit and
$y,z$ indeterminates over $A'$.
Assume there exists an element in $\lambda \in zA_a[[z]]$  which is analytically independent over $A'[[z]]$.
Then the extension $B[[z]]\hookrightarrow B_a[[z]]$ is not TGF
for every   ring $B$ between $A[y]$ and $A'[[y]]$
\end{proposition}
{\noindent\it Proof.}
Let $P$ be the kernel of the $A'_a[[z]]$-algebra morphism
$$\theta:A'_a[[z,y]]\rightarrow A'_a[[z]]$$ given by  $\theta(y)=\lambda$.
Note that $0\neq \lambda-y\in P\cap B_a[[z]]$.
Also $$(P\cap B_a[[z]])\cap B[[z]]\subseteq P\cap {A'}[[z,y]]=0$$ because
$\lambda$ is analytically independent over $A'[[z]]$.
%
$\bullet$

\begin{corollary}\label{8}
Let $D$ be a domain and $x,y,z$ indeterminates over $D$.
Then the extension $B[[z]]\hookrightarrow B_x[[z]]$ is not TGF
for every  ring $B$ between $D[x,y]$ and $D[[x,y]]$.
\end{corollary}
{\noindent\it Proof.}
It suffices to show that there exists an element of $zD[x]_x[[z]]$
which is analytically independent over $D[[x,z]]$, because, after that, we apply
Proposition \ref{3}  for $A=D[x]$,  $A'=D[[x]]$ and $a=x$.

Let $K$ be the quotient field of $D$.
Let $\sigma(z)\in zD[[z]]$ be an element which is algebraically independent over $K(z)$
(e.g., $\sigma(z)=\sum_{n\geq 1} z^{n!}$, cf. \cite[page 277]{B}).
By \cite[Theorem 2]{AM},  $x\sigma(z)$  is analytically independent over $K[[x,xz]]$.
So $x\sigma(z/x)\in zD[x]_x[[z]]$ is analytically independent over $D[[x,z]]$.
$\bullet$\\[2mm]

The next corollary answers \cite[Question 4.8]{HRW}.

\begin{corollary}\label{3c}
Let $K$ be a  field  and $x,y,z$ indeterminates over $K$.
Then the extension $K[x,y][[z]]\hookrightarrow K[x,y]_x[[z]]$ is not TGF.
\end{corollary}
{\noindent\it Proof.}
We apply Corollary \ref{8} for $D=K$ and $B=K[x,y]$.
$\bullet$\\[2mm]

Recall that the {\em generic fiber} of an extension of integral domains $A\subseteq B$ is
the set of prime ideals of $B$ lying over zero in $A$.

\begin{remark}\label{11}{\em
Let $K$ be a  field  and $x,z,y_1,...,y_n$ indeterminates over $K$ with $n\geq 1$.
The arguments employed in the proofs of Proposition \ref{3} and Corollary \ref{8}
can be used to show that the generic fiber of the extension
$$K[x,y_1,...,y_n][[z]]\hookrightarrow K[x,y_1,...,y_n]_x[[z]]$$
has dimension $\geq n$.

Indeed, let
$\sigma_1(z),...,\sigma_n(z)\in zK[[z]]$ be   algebraically independent elements  over $K(z)$.
By \cite[Theorem 2]{AM},  the elements $\lambda_j=x\sigma_j(z/x)$, $1\leq j\leq n$,
are analytically independent over $K[[x,z]]$.
Clearly, $K[[x]]_x$ is the quotient field $K((x))$ of $K[[x]]$.
For $1\leq j\leq n$,
let $P_j$ be the kernel of the $K((x))[[z]]$-algebra morphism
$$\theta:K((x))[[z,y_1,...,y_n]]\rightarrow K((x))[[z,y_{j+1},...,y_n]]$$
given by  $\theta(y_i)=\lambda_i$ for $1\leq i\leq j$ and
$\theta(y_i)=y_i$ for $j+1\leq i\leq n$.
Then $y_j-\lambda_j\in P_j\setminus P_{j-1}$ (in fact,
it is easy to see that $P_j$ is the ideal generated by $y_1-\lambda_1$,...,$y_j-\lambda_j$).
Set $Q_j=P_j\cap K[x,y_1,...,y_n]_x[[z]]$, $1\leq j\leq n$.
Then
$$0=Q_0\subset Q_1\subset Q_2\subset\cdots \subset Q_n$$
the inclusions being proper because $y_j-\lambda_j\in Q_j\setminus Q_{j-1}$.
Since $\lambda_1$,...,$\lambda_n$
are analytically independent over $K[[x,z]]$, we get
$$0=P_n\cap K[[x,z,y_1,...,y_n]]\supseteq Q_n\cap K[x,y_1,...,y_n][[z]].\ $$

Note that the elements $\sigma_1(z),...,\sigma_n(z)$ above  can be chosen in $z\Omega[[z]]$, where
$\Omega$ is the prime subfield of $K$.
Indeed, if $\sigma_1(z),...,\sigma_n(z)\in z\Omega[[z]]$ are
algebraically independent  over $\Omega(z)$,
then, by base extension, they are also algebraically independent over $K(z)$,
because the canonical morphism $\Omega[[z]]\otimes_\Omega K\rightarrow K[[z]]$ is injective
(since every $\Omega$-vector space basis of $K$ is linearly independent over $\Omega[[z]]$).
$\bullet$
}\end{remark}

The next corollary was suggested by the proof of \cite[Proposition 4.9]{HRW}.

\begin{corollary}\label{10}
Let $D$ be a domain, $x,y,z$  indeterminates over $D$
and $a\in D$ a nonzero nonunit such that $\cap_{n\geq 1} a^nD=0$.
Then the extension $B[[z]]\hookrightarrow B_a[[z]]$ is not TGF for every
ring $B$ between $D[x,y]$ and $D[[x,y]]$.
\end{corollary}
{\noindent\it Proof.}
It suffices to show that there exists an element of $zD[x]_a[[z]]$
which is analytically independent over $D[[x,z]]$, because, after that, we apply
Proposition \ref{3}  for $A=D[x]$ and  $A'=D[[x]]$.

By \cite[Theorem 2.1]{S}, there exist $\sigma \in zD[[z/a]]$ which is algebraically independent over
$D[[z]]$.
Now we use the proof of \cite[Proposition 4.9]{HRW}   to show that $\lambda=\sigma x$
is analytically independent over $D[[x,z]]$.
Indeed, let $f(y)\in D[[x,z]][[y]]$, such that $f(\lambda)=0$.
Writing
$f(y)=\sum_{l=0}^\infty \sum_{i+j=l} d_{ij}(z)x^iy^j$ with $d_{ij}(z)\in D[[z]]$,
we get
$$0=f(\lambda)=\sum_{l=0}^\infty \sum_{i+j=l} d_{ij}(z)x^i(\sigma x)^j=
\sum_{l=0}^\infty x^l\sum_{i+j=l} d_{ij}(z)\sigma^j.$$
Hence $\sum_{i+j=l} d_{ij}(z)\sigma^j=0$ for each $l$.
As $\sigma$ is algebraically independent  over $D[[z]]$, each $d_{ij}(z)=0$. Thus $f(y)=0$.
$\bullet$\\[2mm]

The following proposition is the main result of this note.
It is a generalization  of Remark \ref{11}.
\begin{proposition}\label{4}
Let $B$ be a Noetherian domain  containing a field,
$a\in B$  a nonzero nonunit and  $z$  an indeterminate over $B$.
Then the generic fiber of the extension $B[[z]]\hookrightarrow B_a[[z]]$ has
dimension {\em $\geq \mbox{dim}(B/aB)$.}
So, $B[[z]]\hookrightarrow B_a[[z]]$ is not TGF
when {\em $\mbox{dim}(B/aB) \geq 1$.}
\end{proposition}
{\noindent\it Proof.}
Clearly, we may assume that $\mbox{dim}(B/aB)\geq 1$.
Let $1\leq n\leq \mbox{dim}(B/aB)$ and let $M$ be a  prime ideal  of $B$ of height $n+1$ containing $a$.
Set $C=B_M$ and let $\widehat{C}$ be the completion of $C$.
By \cite[Theorem 60]{M}, $\widehat{C}$ contains a coefficient field  $K$.
Pick $b_1,...,b_n\in B$ such that $a,b_1,...,b_n$  is a system of parameters  of $C$.
Then $a,b_1,...,b_n$ is also a system of parameters of $\widehat{C}$.
By the proof of  Cohen's structure theorem for complete local rings \cite[Corollary 2, page 212]{M}
(see also  \cite[Theorem 4.3]{BRD}),
the $K$-algebra  morphism
$$\theta:K[[x,y_1,...,y_n]]\rightarrow \widehat{C}$$ given by
$\theta(x)=a$ and $\theta(y_j)=b_j$, $1\leq j\leq n$, is injective and finite.
We may assume that $\theta$ is the inclusion map. So $x=a$ and $y_j=b_j$, $1\leq j\leq n$.

Let $\Omega$ be the prime subfield of $K$. By Remark \ref{11},
there exist the elements  $\lambda_1,...,\lambda_n\in z\Omega[x]_x[[z]]$
which are analytically independent over $K[[x,z]]$
 and
 a chain of prime ideals of $K((x))[[z,y_1,...,y_n]]$
$$0=P_0\subset P_1\subset P_2\subset\cdots \subset P_n$$
such that $y_j-\lambda_j\in P_j\setminus P_{j-1}$ for  $1\leq j\leq n$ and
$P_n\cap K[[x,z,y_1,...,y_n]]=0$.
Set $D=K[[x,y_1,...,y_n]]$ and
$P'_j=P_j\cap D_x[[z]]$, $1\leq j\leq n$.
Then
$$0=P'_0\subseteq P'_1\subseteq P'_2\subseteq\cdots \subseteq P'_n.$$
%
Since $\widehat{C}_x[[z]]$ is a finite $D_x[[z]]$-module,
there exist a chain of prime ideals of $\widehat{C}_x[[z]]$
$$0\subseteq Q_1\subseteq Q_2\subseteq\cdots \subseteq Q_n$$
such that $P'_j=Q_j\cap D_x[[z]]$ for $1\leq j\leq n$, cf. \cite[Theorem 5]{M}.
Set $T_n:=Q_n\cap \widehat{C}[[z]]$.
Then
$$T_n\cap D[[z]]=
P'_n\cap D[[z]]=P_n\cap D[[z]]=0.$$
As $\widehat{C}[[z]]$ is a finite $D[[z]]$-module, it follows that
$T_n$ is a minimal prime ideal of $\widehat{C}[[z]]$, cf.
 \cite[Theorem 5]{M}. Since $\widehat{C}[[z]]$ is a flat ${B}[[z]]$-module,
we get
$T_n\cap B[[z]]=0$, cf. \cite[Theorem 4]{M}.
Set $N_j=Q_j\cap {B}_x[[z]]$, $1\leq j\leq n$.
Note that $N_n\cap B[[z]]\subseteq Q_n\cap \widehat{C}[[z]]=0$, so
each $N_j$ is in the generic fiber of $B[[z]]\hookrightarrow B_x[[z]]$.
The inclusions
$$0=N_0\subset N_1\subset N_2\subset\cdots \subset N_n$$
are proper  because $y_j-\lambda_j\in N_j\setminus N_{j-1}$ for $1\leq j\leq n$.
Thus the generic fiber of $B[[z]]\hookrightarrow B_x[[z]]$ has dimension $\geq n$.
$\bullet$

\begin{remark}{\em
The assertion of Proposition \ref{4} does not hold for non-Noetherian domains.
For example, let $V$ be a rank-two valuation domain containing a field and let $a$ be
a nonzero element of the height-one prime ideal of $V$.
Then $V_a$ is the quotient field of $V$, so the extension
$V[[z]]\hookrightarrow V_a[[z]]$ is TGF, but $\mbox{dim}(V/aV)=1$.
}\end{remark}

\end{document}